\documentclass[11pt]{article}
\usepackage{epic,latexsym}
\usepackage{amssymb,amsmath,amsthm}
\usepackage{color}
\usepackage{tikz}

\newtheorem{thm}{Theorem}
\newtheorem{lem}[thm]{Lemma}

\newtheorem{ob}[thm]{Observation}
\newtheorem{prop}[thm]{Proposition}
\newtheorem{conj}{Conjecture}

\def \bsk {\bigskip}

\def \nmr {\begin{enumerate}}
\def \enmr {\end{enumerate}}

\def \tmz {\begin{itemize}}
\def \etmz {\end{itemize}}

\newcommand{\Proof}{\noindent\textbf{Proof. }}
\newcommand{\cH}{{\cal H}}
\newcommand{\cP}{{\cal P}}

\newcommand{\tg}{\gamma_{{\rm tg}}}
\newcommand{\gtd}{\gamma_{{\rm tg}}}

\begin{document}

\title{On the game total domination number\thanks{Research supported by the
       National Research, Development and Innovation
Office -- NKFIH under the grant SNN 116095.}
}
\author{Csilla Bujt\'{a}s
\\ \\
Faculty of Information Technology \\
University of Pannonia, Veszpr\'em, Hungary\\
\small \tt Email: bujtas@dcs.uni-pannon.hu,
   \\
}
\date{}
\maketitle

\begin{abstract}
The total domination game is a two-person
competitive optimization game, where the players, Dominator and
Staller, alternately select vertices of an isolate-free graph $G$.
Each vertex chosen must strictly increase the number of vertices
totally dominated. This process eventually produces a total
dominating set of $G$. Dominator wishes to minimize the number
of vertices chosen in the game, while Staller wishes to maximize it.
The game total domination number of $G$, $\gamma_{{\rm tg}}(G)$, is
the number of vertices chosen when Dominator starts the game and
both players play optimally.

Recently, Henning,  Klav\v{z}ar, and Rall proved that $\gamma_{{\rm
tg}}(G) \le \frac{4}{5}n$ holds for every graph $G$ which is given
on $n$ vertices such that every component of it is of order at least
$3$; they also conjectured that the sharp upper bound would be
$\frac{3}{4}n$.
 Here, we prove that $\gamma_{{\rm tg}}(G) \le
 \frac{11}{14}n$ holds for every $G$ which contains no isolated
 vertices or isolated edges. \\

\end{abstract}


\section{Introduction}

Total domination game is a two-person competitive optimization game based on
 the notion of total domination.
We study the corresponding graph invariant $\tg(G)$, called game
total domination number. Our main contribution is a general upper
bound  $\frac{11}{14}n$ on $\tg(G)$ that holds for every graph $G$
of order $n$ not containing isolated vertices or isolated edges. In
the proof, we will consider the so-called `open neighborhood
hypergraph' $\cH$ instead of $G$, and assign weights to the vertices
and edges of $\cH$. Then, we analyze a greedy strategy of the `fast'
player, called Dominator.

\subsection{Basic terminology}

For a graph  $G$ and a vertex $v\in V(G)$, the open neighborhood of
$v$ is $N_G(v)= \{u: uv \in E(G)\}$, and its closed neighborhood is
$N_G[v]=N_G(v) \cup \{v\}$. If $S\subseteq V(G)$, $N_G(S)=\bigcup_{v
\in S} N_G(v)$ and $N_G[S]=\bigcup_{v \in S} N_G[v]$. We say that a
vertex $v$ \emph{totally dominates} $u$ if  $u \in N_G(v)$, while
$v$ \emph{dominates} $u$ if  $u\in N_G[v]$. A set $D$ of vertices is
a \emph{total dominating set} and a \emph{dominating set} in $G$ if
  $N_G[D]=V(G)$ and $N_G(D)=V(G)$ holds respectively. Equivalently, $D$ is a total dominating set if
  each vertex  has a neighbor in $D$, and $D$ is a
  dominating set if
 each vertex  which is not in $D$ has a neighbor in $D$.
 The invariant \emph{total domination number} $\gamma_t(G)$ and  \emph{domination number} $\gamma(G)$
 is the minimum size of a total dominating set and that of a dominating
 set in $G$, respectively.

The notion of \emph{total domination game} was introduced recently
by Henning,  Klav\v{z}ar, and   Rall \cite{HKR-2015}. It is played
on an isolate-free  graph $G$ by two players, namely Dominator and
Staller, who alternately select vertices of $G$.
 A move (a selection) is legal if the chosen vertex totally dominates at least one vertex
 which is not totally dominated by the set of vertices previously selected.
  The game is over when the set $D$ of chosen
 vertices becomes a total dominating set in $G$.
  Dominator wishes to finish the game as soon as possible, while Staller
 wishes to delay the end of the game. The \emph{game total domination number}, $\tg(G)$, of $G$ is the number of vertices chosen
    when   Dominator starts the game and both  players play
  optimally.\footnote{ Remark that the total domination game is an analogous version of
  the domination game, introduced by
Bre\v{s}ar, Klav\v{z}ar, and Rall in 2010 \cite{BKR-2010},
   where the choice of  $v$ is legal if it dominates
 at least one new vertex; that is, if  $N[v]\setminus N[D]\neq
 \emptyset$.
 The corresponding invariant is the game domination number,
 $\gamma_g(G)$.
 For the exact definitions and results on the domination game see
 \cite{BKR-2010, KWZ-2013, BDKK-2014}.}

    A \emph{hypergraph} $\cH$ is a set (multi)system over the vertex set $V(\cH)$. The edge set $E(\cH)$ of  $\cH$ contains nonempty subsets of
    $V(\cH)$. An edge $e\in E(\cH)$ is a \emph{$k$-edge} if $|e|=k$. $\cH$ is a \emph{linear hypergraph}, if for any two different edges
    $e_1, e_2 \in E(\cH)$,
        $|e_1 \cap e_2|\le 1$ holds. In particular,  there are no
        multiple edges of size greater than 1 in a linear
        hypergraph.
        The degree $d_\cH(v)$
    of a vertex $v\in V(\cH)$ is the number of edges incident to $v$, and the maximum degree $\Delta(\cH)$ equals
    $\max\{d_\cH(v): v\in V(\cH)\}$. A \emph{vertex cover} (also called
    transversal) in $\cH$ is a set $T$ of vertices which contains at
    least one vertex from each edge.
    Remark that unlike the usual terminology, here we allow also the presence of multiple edges and $1$-edges in hypergraphs.

    \emph{Transversal game} on hypergraphs was introduced recently in \cite{BHT-2016} and further studied in \cite{BHT-2017}.
    Its definition is analogous to that of total domination game. Two players, namely Dominator\footnote{The `fast' player is called Edge-hitter in \cite{BHT-2016} and \cite{BHT-2017}. To have the two players with the same names in the total domination and in the transversal game, we prefer to call him Dominator instead of Edge-hitter, here.}
    and Staller, alternately choose vertices of a hypergraph $\cH$. A move is legal if the vertex chosen covers at least one edge which has not been covered in the game so far. The game is over when every edge of $\cH$ is covered. Dominator wishes to end the game as soon as possible, while Staller wishes to delay the end of the game. Assuming that Dominator starts the game on $\cH$, and also that both players play optimally, the length of the game is uniquelly determined. It is called the \emph{game transversal number} of $\cH$ and denoted by $\tau_g(\cH)$.

 Given an isolate-free graph $G$, its \emph{open neighborhood hypergraph} $ONH(G)$ is the hypergraph with vertex set $V(G)$ and edge set
$$E(ONH(G))=\{N_G(v): v \in V(G)\}.$$
It is easy to see (and also observed earlier) that a vertex set $T$
is a total dominating set in $G$ if and only if it is a vertex cover
in $ONH(G)$. Similarly, a sequence of moves defines a legal total
domination game on $G$, if and only if it is a legal transversal
game on $ONH(G)$. Consequently, $\tg(G)=\tau_g(ONH(G))$ holds for
every isolate-free graph $G$.

\subsection{Results}

In the introductory paper \cite{HKR-2015},  among other basic
results, the sharp bounds $\gamma(G)\le \tg(G) \le 3\gamma(G)-2$ are
proved. The exact value of $\gtd(G)$
 for paths and cycles were established in \cite{DH}.
 Our present subject is strongly connected to the following `$\frac{3}{4}$-Game Total Domination Conjecture',
 posed by Henning,  Klav\v{z}ar, and   Rall  in \cite{HKR-2017}.

\begin{conj}
\label{conj} If $G$ is a graph on $n$ vertices in which every
component contains at least three vertices, then $\gtd(G) \le
\frac{3}{4}n$.
\end{conj}
Note that the restriction given on the size of the components is
necessary, because otherwise the upper bound on $\tg(G)$ could not
be better than $n$.
 We also remark that if the conjecture is true then it is sharp. Tight examples
given in \cite{HKR-2017} are the graphs each component of which is a
path of length 4 or 8.

 The following results related to Conjecture~\ref{conj}
have been proved so far. In each of them, it is assumed that $G$ is
a graph of order $n$ in which every component contains at least
three vertices.
  \tmz
    \item  $\gtd(G) \le  \frac{4}{5}\; n$. \cite{HKR-2017}
    \item  If $\delta(G) \ge 2$, then $\gamma_{tg}(G) < \frac{8}{11}\; n < \frac{3}{4}\;
    n$.  \cite{BHT-2016}
    \item  If $\deg(u)+\deg(v) \ge 4$ for every edge $uv \in E(G)$, and no two
  degree-1 vertices are at distance 4, then $\gamma_{tg}(G) \le
  \frac{3}{4}\;n$.  \cite{HR-2016}

    \etmz

In this paper our main contribution is a new general upper bound on
the game total domination number that improves the earlier bound
$4n/5$.

\begin{thm} \label{thm-1}
 If $G$ is a graph of order $n$ in which every
component contains at least three vertices, then $$\gtd(G) \le
\frac{11}{14}\;n.$$
\end{thm}

This theorem will be proved in Section~\ref{sec-2}. In
Section~\ref{sec-3}, we make some concluding remarks on the
Staller-start version of the total domination game.

\section{Proof of the upper bound $11n/14$}
\label{sec-2}

In this section we prove our main result, namely
Theorem~\ref{thm-1}. Given a graph $G$ which does not contain
isolated vertices and isolated edges, we construct its open
neighborhood hypergraph $\cH_0=ONH(G)$. Then, the total domination
game on $G$ will be represented by the transversal game on the
hypergraph $ONH(G)$ where the same sequence of vertices is played.
Since our aim is to give a general upper bound on $G$, degree-1
vertices in $G$ and the corresponding edges of size
one in $ONH(G)$ are not excluded. 
Throughout the proof, we denote by $j^*$ the number of turns in the
game. Let $m_k$ be the vertex chosen in the $k$th turn  ($1\le k \le
j^*$).
 We set $D_0=\emptyset$ and define $D_i=\{m_k: 1\le k \le i\}$ for $1\le i \le
 j^*$.

\subsection{Residual hypergraph and special vertices}

During the transversal game, the edges which are already covered and
the vertices which are not incident with any uncovered edges do not
influence the continuation of the game. Hence, we delete them and
obtain the residual hypergraph $\cH_i$. It is defined formally as
$$E(\cH_i)=\{e\in E(\cH_0): e\cap D_i = \emptyset \} \quad \mbox{and} \quad V(\cH_i)=\bigcup_{e \in E(\cH_i)} e,$$
where $\cH_0$ denotes $ONH(G)$. Note that $\cH_{j^*}$ is the empty
hypergraph.

Roughly, we would  like to say that a vertex $v$ and the corresponding edge $e_v= N_G(v)$ is special in
 $ONH(G)=\cH_0$ if   $d_G(v)=d_{\cH_0}(v)=1$ and consequently, $e_v$ is a 1-edge in the hypergraph.
 But we do not need more than one special vertex inside any edge of $\cH_0$.
 So, the definition will be the following. Consider all the edges of $\cH_0$
 that contains at least one degree-1 vertex and (arbitrarily) fix exactly one degree-1
 vertex from each such edge. These vertices will be referred to as  \emph{special vertices}
 and the set of the special vertices will be denoted by $S$. If $v\in S$,
 then the corresponding edge $e_v$ is called \emph{special edge}. The definitions imply the following simple statements.
\begin{ob}
\label{obs} Let $G$ be a graph which contains no isolated vertices
and isolated edges,  $\cH_0$ be its open neighborhood hypergraph,
 and $S$ be a  fixed set of special vertices in $\cH_0$. Let $\cH_i$ be the residual hypergraph obtained in a transversal game on $\cH_0$ ($0\le i $).
\tmz
  \item[$(i)$] The number of special vertices equals the number of special edges in $\cH_0$.
    \item[$(ii)$] Any edge of $\cH_i$ contains at most one special vertex.
    \item[$(iii)$] Any vertex in $\cH_i$ is incident with at most one special edge.
    \item[$(iv)$] No special edge contains a special vertex in $\cH_i$.
\etmz
\end{ob}
\Proof The definitions immediately imply that the statements
$(i)$-$(iii)$ are valid for $\cH_0$. Moreover, if $\cH_0$ satisfies
$(ii)$ and $(iii)$, these remain valid for every later residual
hypergraph $\cH_i$.
  Concerning $(iv)$, we observe that a special edge containing a special
  vertex in $\cH_0$ would correspond to a degree-1 vertex in $G$ the neighbor of which is also of degree 1.
  This contradicts the exclusion of $P_2$-components from $G$.
  Since $\cH_0$ satisfies $(iv)$,  every  residual hypergraph $\cH_i$ ($i \ge 1$) satisfies it as well.
    \qed

We emphasize that an edge or a vertex is special in $\cH_i$, if it
was special  in $\cH_0$ and it is still present in the residual
hypergraph $\cH_i$.

\subsection{Weights and phases}

In a residual hypergraph $\cH_i$, a component will be called
\emph{Type-X component}, if it corresponds to an isolated edge which
contains at least two non-special vertices. The number of Type-X
components in $\cH_i$ is denoted by $x_i$.   Moreover, $n_i^h$ and
$n_i^s$ denote the number of non-special and special vertices
present in $\cH_i$, while $e_i^h$ and $e_i^s$ stand for the number
of non-special and special edges present in $\cH_i$, respectively.
 We define the following function on the residual hypergraphs
$$f(\cH_i)=13 n_i^h+ 7n_i^s + 9e_i^h + 15e_i^s-7x_i.$$
In an equivalent formulation, we may say that the following weights
$f(v)$ and $f(e)$ are assigned to every vertex $v$ and edge $e$.
\begin{center}
\begin{tabular}{|c|c|c|}  \hline
 &   Non-special & Special  \\
\hline
 Vertex $v$ &  $f(v)= 13$ &  $f(v)$= 7\\ \hline
 Edge $e$ &  $f(e)=9$ &  $f(e)=15$\\ \hline
\end{tabular}
\end{center}
Then, the weight of the residual hypergraph is
 $$f(\cH_i)= \sum_{v\in V(\cH_i)}f(v)+\sum_{e \in
 E(\cH_i)}f(e)-7x_i.$$

Since every $v\in V(\cH_0)$ and  the corresponding edge $e_v$ satisfies $f(v)+f(e_v)=22$, we have
 $f(\cH_0) \le 22\;n$
and $f(\cH_{j^*})=0$. In the $i$th turn of the game, $1\le i \le
j^*$, the \emph{decrease in the weight} is
$d_i=f(\cH_{i-1})-f(\cH_i)$. We will suppose that Dominator follows
a greedy strategy in the transversal game; that is, for every odd
$i$, he plays a vertex in $\cH_{i-1}$ which results in the possible
maximum decrease  $d_i$. Our aim is to prove that, under this greedy
strategy,
  $$\frac{\sum_{i=1}^{j^*} d_i}{j^*}\ge 28$$
    is always valid for the average decrease in a turn, independently of Staller's strategy.

To analyze the game, we split it into four phases. Let $[j^*]$
denote $\{1,\dots, j^*\}$ and define the following sets
\[
\begin{array}{lcl}
\cP^1 & = & \{i \in [j^*]: \forall \ell \; ((\ell \mbox{ is odd and } \ell \le i) \rightarrow \; d_\ell \ge 40)\}, \\
\cP^2 & = & \{i \in [j^*]: \forall \ell \; ((\ell \mbox{ is odd and } \ell \le i) \rightarrow \; d_\ell \ge 38)\} \setminus \cP^1,\\
\cP^3 & = & \{i \in [j^*]: \Delta(\cH_{i-1})\ge 2\} \setminus (\cP^1 \cup \cP^2),\\
\cP^4 & = & \{i \in [j^*]: \Delta(\cH_{i-1})=1\} \setminus (\cP^1
\cup \cP^2).
\end{array}
\]
 By definition,  each $\cP^k$ (if not empty) contains consecutive
integers. Moreover, $\{\cP^1,\cP^2, \cP^3, \cP^4\}$ gives a
partition of $[j^*]$. We say that the $i$th turn of the game belongs
to Phase $k$ if $i \in \cP^k$. To simplify the later formulas, we
also define $a_k=\min (\cP^k)$ and $b_k=\max (\cP^k)$ if $\cP^k$ is
not empty. If $\cP^k$ is empty, we define $b_k=b_{k-1}$ ($1\le k \le
4$) artificially that can be done recursively if we set $b_0=0$.

\subsection{Phase 1}

At the beginning of this subsection, we prove two general lemmas
that remain valid throughout the game. The first of them gives a
lower bound on the decrease of the weights, if a vertex from an
isolated edge is played.

\begin{lem}
\label{lem-iso}
  If $i \in  [j^*]$,  $m_i=v$ and $v$ belongs to an isolated edge $e$ in
  $\cH_{i-1}$, then $d_i \ge 28$. In particular, $d_i \ge 28$, if
  $v$ is from a Type-X component.
\end{lem}
\Proof Since $e$ is an isolated edge, after the move $m_i=v$, the
edge $e$ and all vertices from it will be deleted. If $e$ is a
1-edge, $e$ is special and, by Observation~\ref{obs}$(iv)$, $v$ is
not special. Thus, $d_i \ge 15+13=28$. If $e$ is a 2-edge and
contains a special vertex, the other vertex is not special and
hence, $d_i \ge 9+7+13=29$. In the remaining cases $e$ contains at
least two non-special vertices; that is, $e$ is from a Type-X
component. This means $x_i=x_{i-1}-1$, and the decrease is $d_i \ge
9 + 2\cdot 13-7 =28.$ \qed

 Next, we prove a lower bound on the
decrease $d_i$ in the weight. It is true  regardless of that the
next player is Dominator or Staller.
\begin{lem}
\label{St-16}
  For every $i\in [j^*]$, $d_i \ge 16.$
\end{lem}
\Proof In the $i$th turn, at least one new edge is covered and
deleted from the hypergraph $\cH_{i-1}$; and at least one vertex
(the one which was played) is deleted. Hence, if $x_i \ge x_{i-1}$,
then $d_i \ge 9+7=16$. If $x_i < x_{i-1}$, then a vertex from a
Type-X component was played. By Lemma~\ref{lem-iso}, we have $d_i
\ge 28$ that completes the proof.
  \qed

  Now we are ready to prove that the average decrease in a turn is at least
$28$ in Phase~1.
 \begin{lem}
 \label{lem-ph1}
  If  $\cP^1\neq \emptyset$,
    $$\frac{\sum_{i=1}^{b_1} d_i}{b_1}\ge 28.$$
\end{lem}
\Proof If $i$ is odd and $1\le i \le b_1$, the definition of $\cP^1$
ensures that $d_i \ge 40$. By Lemma~\ref{St-16}, we have
$d_i+d_{i+1}\ge 40+16= 2\cdot 28$ that implies the statement. Remark
that if the game finishes with Dominator's turn in  Phase 1 (i.e.,
$b_1$ is odd and equal to $j^*$), then the last decrease $d_{b_1}$
is at least $40$, and $\sum_{i=1}^{b_1} d_i \ge \frac{b_1-1}{2}\cdot
56 +40 > 28 \;b_1$. Hence, the lemma is valid for this special case
as well. \qed

We may prove some properties which are true for each residual
hypergraph after the end of Phase 1.
\begin{lem}
\label{lem-end1} For every $i \ge b_1$, the residual hypergraph
$\cH_i$ satisfies the following properties.
  \tmz
    \item[$(i)$] $\Delta(\cH_i) \le 2$.
    \item[$(ii)$]  $\cH_i$ is a linear hypergraph.
    \item[$(iii)$] If  $v$ is special vertex and $u$ is a
    neighbor of $v$ in $\cH_i$, then $u$ is not contained in any
    special edges.
     \item[$(iv)$] If $d_{\cH_i}(v)=2$, then $|S\cap
     N_{\cH_i}(v)|\le 1$.
    \etmz
\end{lem}
\Proof $(i)$ Assume for a contradiction that there exists a vertex
$v\in V(\cH_{b_1})$ which is incident with at least three edges.
Then Dominator may play $v$ in the next turn and $\cH_{b_1+1}$ is
obtained from $\cH_{b_1}$ by deleting $v$ which is not special (and
maybe, some further vertices),  and at least three edges. Clearly,
$x_{b_1+1}\ge x_{b_1}$. Thus, we have $d_{b_1+1}\ge 13+3 \cdot 9=40$
that would imply  $b_1+1 \in \cP_1$ that is a contradiction. Hence,
$\Delta(\cH_i) \le 2$ holds for $i=b_1$ and for every larger index.

$(ii)$ Now, assume that there exist two edges, say  $e_1$ and $e_2$,
in $\cH_{b_1}$ such that $|e_1 \cap e_2|\ge 2$. By $(i)$, every
vertex from $e_1\cap e_2$ is of degree 2. Playing a common vertex
$v$ of $e_1$ and $e_2$, at least two non-special vertices and two
edges will be deleted. Since the number of Type-X components is not
decreased, $d_{b_1+1} \ge 2\cdot 13+ 2\cdot 9=44$ that is a
contradiction, again. This proves $(ii)$ for $i = b_1$ and implies
the linearity for every later residual hypergraph.

$(iii)$ Assume for a contradiction that $v\in S$, $u\in
N_{\cH_{b_1}}(v)$ and $e=\{u\}$ is a special edge. By
Observation~\ref{obs}$(ii)$, $u$ is not a special vertex. If
Dominator plays $u$, then $u$, $v$, $e$, and the edge incident with
$v$
 will be deleted. Since $x_{b_1+1}\ge x_{b_1}$,
  we have $d_{b_1+1} \ge 13+7+15+9=44
>40$, a contradiction. As new special vertices and edges cannot
arise during the game, $(iii)$ holds for every $i \ge b_1$.

$(iv)$ Consider first $\cH_{b_1}$. Suppose that a vertex $v$
 is incident with two edges $e_1$, $e_2$ and that $N_{\cH_{b_1}}(v)$ contains
 two special vertices $u_1$ and $u_2$. If Dominator plays $v$ in the
 next turn, the vertices $v$, $u_1$ ,$u_2$ and the edges $e_1$,
 $e_2$ will be deleted from $\cH_{b_1}$ and $x_{b_1+1}\ge x_{b_1}$.
 This would yield $d_{b_1} \ge 13 +2\cdot 7+2\cdot 9=45 >40$ that is
 a contradiction. This proves the statement for $i=b_1$ from which
 $(iv)$ follows for every $i\ge b_1$.
 \qed

\subsection{Phase 2}

In Phase 2, every residual hypergraph satisfies the properties
$(i)$-$(iv)$ from Lemma~\ref{lem-end1}, and $d_i \ge 38$ for every
odd $i$. We will prove that the average decrease is at least 28 over
the turns in Phase 2.

\begin{lem}
\label{lem-ph2}
  If  $\cP^2\neq \emptyset$,
    $$\frac{\sum_{i=a_2}^{b_2} d_i}{|\cP^2|}\ge 28.$$
\end{lem}

 \Proof Suppose that Staller plays a vertex $v$ in
the $i$th turn, $a_2<i \le b_2$. By definition of $\cP^2$, $d_{i-1}
\ge 38$ holds. We have the following cases concerning the move
$m_i=v$.
  \tmz
  \item  If $v$ is from an isolated edge, then, by Lemma~\ref{lem-iso}, $d_i \ge
  28$. This gives $d_{i-1}+d_i \ge 38+28 >2\cdot 28.$
  \item If $v$ is a non-special vertex and it is not from an
  isolated edge (i.e., $x_i \ge x_{i-1}$), then $d_i \ge 13+9=22$ and again, we have $d_{i-1}+d_i \ge 38+22 >2\cdot
  28$.
  \item In the third case, $v$ is a special vertex   and it is
  contained in a non-isolated edge in $\cH_i$. Then, $v$ has a degree-2
  neighbor, say $u$. Let $e$ be the edge containing both $v$ and
  $u$, and let $e'$ be the other edge incident with $u$. By
  Lemma~\ref{lem-end1}$(iv)$, $e'$ is not special. This implies $|e'|\ge 2$. Moreover, by
  Lemma~\ref{lem-end1}$(iii)$, $e'$ does not contain any special
  edges. Consequently, $e'$ contains at least two non-special
  vertices.
  If $e'$ becomes isolated in $\cH_{i+1}$, then it will be of a Type-X edge. Therefore, Staller's move
  results in a decrease of $d_i \ge 7+9+7=23$, since $v$ and $e$ are
  deleted and $x_{i+1}\ge x_i+1$. Hence, we have, $d_{i-1}+d_i \ge 38+23 >2\cdot
  28$.
  In the other case, $e'$ is not isolated and contains a degree-2 vertex $u'$ in
  $\cH_{i+1}$. Then, Dominator may choose $u'$ in the $(i+1)$st
  turn. This means that two non-special vertices, namely $u$,
  $u'$, and two edges are deleted from the residual hypergraph. Hence, $d_{i+1}\ge 2\cdot 13 +
  2\cdot 9 =44$. If the game is finished with the $(i+1)$st turn,
  then $d_{i-1}+d_i+d_{i+1}\ge 38+16+44>3\cdot 28$. In the other case,  $\cH_{i+1}$
  is not empty, and we have $d_{i-1}+d_i+d_{i+1}+d_{i+2}\ge 38+16+44+16>4\cdot 28$.
  \qed
  \etmz

  \begin{lem}
\label{lem-end2} For every $i \ge b_2$, every special vertex present
in $\cH_i$ is contained in an isolated edge.
\end{lem}
\Proof Consider $\cH_{b_2}$ and assume that the special vertex $v$
has a degree-2 neighbor $u$. If Dominator plays $u$ in the next
turn, then $v$, $u$ and two edges are deleted, moreover
$x_{b_2+1}\ge x_{b_2}$. This gives $d_{b_2+1}\ge 7+13+2\cdot 9=38$
that contradicts the definition of $\cP^3$. Thus, the lemma is valid
with $i=b_2$ and in turn, it implies the statement for every later
residual hypergraph. \qed

\begin{lem}
\label{lem-22}
  For every $i \ge b_2+1$, $d_i \ge 22$ holds.
\end{lem}
\Proof After the end of Phase 2, in every turn, either an isolated
edge is deleted that gives $d_i \ge 28$ by Lemma~\ref{lem-iso}, or a
vertex $v$ is played which does not belong to an isolated edge. In
the latter case, by Lemma~\ref{lem-end2}, $v$ cannot be special and
the move results in a decrease $d_i\ge 13+9=22$. \qed

\subsection{Phase 3 and 4}
If Phase $3$ is not empty, it starts with the turn $a_3$. By the
definition of the phases, Dominator's greedy strategy gives $d_{a_3}
<38$, while $d_i \ge 38$ for every odd $i$ smaller than $a_3$.
Moreover, $\Delta(\cH_k) =2$ holds for each $a_3-1 \le k \le b_3-1$.
Therefore, in Phase 3, Dominator can always play a vertex of degree
2 that results in a decrease of at least $13+2\cdot 9= 31$ in the
weight of the residual hypergraph.

\begin{lem}
  \label{lem-ph3}
  If  $\cP^3\neq \emptyset$,
    $$\frac{\sum_{i=a_3}^{b_3} d_i}{|\cP^3|}\ge 28.$$
\end{lem}

  \Proof Consider an odd $i$ with $a_3 \le i \le
b_3$. First suppose that there exists a degree-2 vertex $v$ which
has a degree-1 neighbor $u$ in $\cH_{i-1}$.   Remark that by
Lemma~\ref{lem-end2}, both $v$ and $u$ are non-special vertices.
Then, Dominator may play $v$, and this move results in $d_i \ge
2\cdot 13 + 2\cdot 9 =44$. Since by Lemma~\ref{lem-22}, $d_{i+1}\ge
22$, we have $d_i +d_{i+1}\ge 44+22>2\cdot 28$ in this case.

Second, suppose that every component of $\cH_{i-1}$ which is not an
isolated edge is 2-regular. If Dominator may play a vertex such that
a new isolated edge arises, then $d_i \ge 13+ 2\cdot 9+7=38$ and
$d_i +d_{i+1}\ge 38+22>2\cdot 28$ follows. Also, if there exists a
special edge on the degree-2 vertex $v$, then the choice of $v$
gives $d_i \ge 13+15+9=37$ and in turn, $d_i +d_{i+1}\ge
37+22>2\cdot 28$. In the remaining case, Dominator plays a degree-2
vertex $v$ which has at least two  neighbors, say $u_1$ and $u_2$.
Further, as new isolated edges do not arise, $u_1$ and $u_2$ become
degree-1 vertices in a component the maximum degree of which is two.
This results in $d_i \ge 31$. If Staller's move creates a new
isolated edge, $d_{i+1}\ge 13+9+7=29$. If after Staller's turn both
$u_1$ and $u_2$ are deleted from the residual graph, $d_{i+1}\ge
2\cdot 13+9=35$. In both cases $d_i+d_{i+1}\ge 31+29
>2\cdot 28$. So, it is enough to consider the case when $d_i \ge
31$, $d_{i+1}\ge 22$, and at least one of $u_1$ and $u_2$, say
$u_1$,  is a degree-1 vertex contained in a non-isolated edge of
$\cH_{i+1}$. Thus, $u_1$ has a neighbor $w$ of degree 2 in
$\cH_{i+1}$. If Dominator selects $w$ in the next turn, $u_1$, $w$
and two edges will be deleted. Hence, $d_{i+2}\ge 2\cdot 13+2\cdot
9=44$. If the game is finished with this turn,
$d_i+d_{i+1}+d_{i+2}>3\cdot 28$ holds. If the game continues with
the $(i+3)$rd turn,
$$d_i+d_{i+1}+d_{i+2}+d_{i+3}\ge 31+22+44+22=119 >4\cdot 28$$
follows. This finishes the proof of the lemma.
 \qed

\begin{lem}
\label{lem-ph4}
  If  $\cP^4\neq \emptyset$,
    $$\frac{\sum_{i=a_4}^{b_4} d_i}{|\cP^4|}\ge 28.$$
\end{lem}

 \Proof In Phase 4, by definition, we have only
isolated edges and by Lemma\ref{lem-iso}, $d_i \ge 28$ follows for
every $i\ge a_4$ in the game. This proves the lemma. \qed

\bsk
 \noindent By  Lemma~\ref{lem-ph1},
Lemma~\ref{lem-ph2}, Lemma~\ref{lem-ph3}, and Lemma~\ref{lem-ph4},
we have that
 $$\frac{\sum_{i=1}^{j^*} d_i}{j^*}= \frac{f(\cH_0)-f(\cH_{j^*})}{j^*}\ge 28,$$
 where $j^*$ denotes the length of the game when Dominator follows a
 greedy strategy based on the function $f$, and Staller plays optimally, according to her goal. Consequently,
 $$\tg(G)=\tau_g(\cH_0)\le j^* \le \frac{f(\cH_0)}{28}\le \frac{11}{14}\; n$$
follows, which proves Theorem~\ref{thm-1}.

\section{Concluding remarks}
\label{sec-3}

Analogously to the game total domination number $\tg(G)$ (resp., to the game transversal number $\tau_g(\cH)$), the \emph{Staller-start game total
  domination number}, $\tg'(G)$ (resp., the \emph{Staller-start game transversal number} $\tau_g(\cH)$)is the length of the game if Staller
  starts and both players play optimally. It was proved already in the introductory paper \cite{HKR-2015} that for  any graph $G$,
    $|\tg(G)-\tg'(G)|\le 1$.

    In \cite{HKR-2017}, the authors also posed a conjecture on the Staller-start version of the total domination game.

    \begin{conj}
\label{conj} If $G$ is a graph on $n$ vertices in which every
component contains at least three vertices, then $\gtd'(G) \le
\frac{3n+1}{4}$.
\end{conj}

In the same paper, they proved the upper bound $\frac{4n+2}{5}$. Our
proof given for Theorem~\ref{thm-1} can be easily extended with a
Preliminary Phase which contains the first move $m_0$ taken by
Staller. In this short part, $f(ONH(G))$ is decreased by $d_0\ge 16$
and then, from the next turn, the determination of the phases and
the proof is just the same as it was in Section~\ref{sec-2}.
Consequently, we have
 the upper bound $\frac{22n+12}{28}=\frac{11n+6}{14}$.
\begin{prop}
 If $G$ is a graph on $n$ vertices in which every
component contains at least three vertices, then $\gtd'(G) \le
\frac{11n+6}{14}$.
\end{prop}

\end{document}